\documentclass[12pt]{amsart}
\usepackage[top=1.2in, left=1.2in, right=1.2in, bottom=1.2in]{geometry}
\usepackage[utf8x]{inputenc}
\usepackage{amsmath,amsfonts,amssymb,epsfig}
\usepackage{delarray}
\usepackage{graphicx}
\usepackage{url}
\usepackage{framed}

\newtheorem{theorem}{Theorem}[section]

\newtheorem{cor}[theorem]{Corollary}
\newtheorem{definition}[theorem]{Definition}

\newtheorem{example}[theorem]{Exemple}

\newtheorem{remark}[theorem]{Remark}

\newcommand{\bbR}{\mathbb{R}} %real numbers
\newcommand{\cB}{\mathcal{B}}

\newcommand{\cD}{\mathcal{D}}
\newcommand{\cE}{\mathcal{E}}

\newcommand{\cG}{\mathcal{G}}

\newcommand{\cJ}{\mathcal{J}}
\newcommand{\cL}{\mathcal{L}}

\newcommand{\cO}{\mathcal{O}}

\newcommand{\cV}{\mathcal{V}}

\newcommand{\intprod}{\!\mathbin{\hbox to 6pt
{\vrule height0.4pt width5pt depth0pt \kern-.4pt
\vrule height6pt width0.4pt depth0pt\hss}}\!}

\begin{document}

\title[Integrable systems in planar robotics]{Integrable systems in planar robotics}

\author{Tudor S. Ratiu}
\address{School of Mathematical Sciences, Shanghai Jiao Tong
University, 800 Dongchuan Road, Minhang District, Shanghai, 200240 China,
Section de Math\'ematiques, Universit\'e de Gen\`eve, 2-4 rue du
Li\`evre, Case postale 64, 1211 Gen\`eve 4, and Ecole Polytechnique
F\'ed\'erale de Lausanne, CH-1015 Lausanne, Switzerland.  }
\email{ratiu@sjtu.edu.cn}
\thanks{Partially
supported by the National Natural Science Foundation of China  grant
number 11871334 and by NCCR SwissMAP grant of the Swiss National Science Foundation.}

\author{Nguyen Tien Zung}
\thanks{Partially
supported by the ``High-End Foreign Experts'' cooperation program at
Shanghai Jiao Tong University.}
\address{Institut de Math\'ematiques de Toulouse, UMR5219, 
Universit\'e Toulouse 3}
\email{tienzung@math.univ-toulouse.fr}

\begin{abstract}
The main purpose of this paper is to investigate commuting flows 
and integrable systems on the configuration 
spaces of planar linkages. Our study leads
to the definition of a natural volume form on each configuration
space of planar linkages, the notion of cross  products
of integrable systems,  and also the notion of multi-Nambu integrable systems.
The first integrals of our systems are functions of Bott-Morse type, 
which may be used to study the topology of configuration spaces. \\

\centerline{\textit{Dedicated to Anatoly T. Fomenko on the occasion 
of his 75th birthday}}
\end{abstract}

\maketitle

%\tableofcontents

\section{Introduction}

\begin{figure}[!ht]
\centering
\includegraphics[height = 0.36 \linewidth]{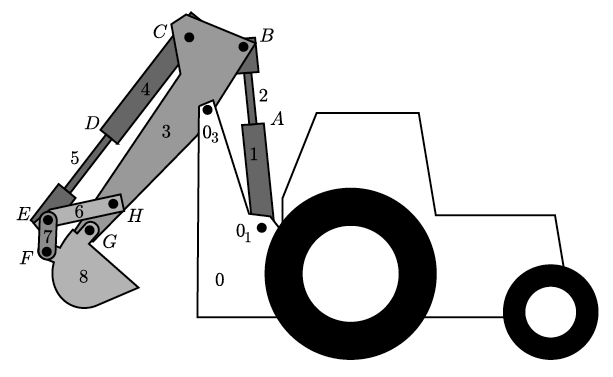} 
\includegraphics[height = 0.36 \linewidth]{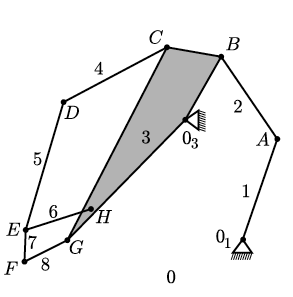} 
\caption{The excavator and its planar linkage representation,
after \cite{Servatius-Assur2010}.}
\label{fig:Excavator}
\end{figure}

This paper is a preliminary investigation of natural commuting flows and
integrable systems on the configuration spaces of linkages in geometry, 
machinery and robotics, especially planar linkages. Our motivation is manifold:

$\bullet$ Firstly,  linkages appear everywhere in mechanical systems 
(see, e.g., \cite{McCarthy-Designs2011,ZFMC-Design}).
Moreover, commuting flows may be viewed as ``independent components'' 
in a mechanical system  and play an important role in robotics and motion 
control: the commutativity makes it easier to compute and control 
the motions. In fact, one may observe many commuting flows in practice. For example,
the excavator in Figure \ref{fig:Excavator} (left), borrowed from Servatius et al. 
\cite{Servatius-Assur2010} with permission,
has two sliding mechanisms, which commute with each
other. This excavator may be represented by
a planar linkage in Figure \ref{fig:Excavator} (right), 
where each sliding mechanism is replaced
by two connecting edges. The configuration space of this linkage admits an 
integrable system having two commuting vector fields. If one allows some of the 
lengths of the edges to be variable, then the configuration space has higher 
dimension, and those adjustable lengths are first integrals of the integrable 
system in question.

$\bullet$ Secondly, for 3-dimensional (3D) polygonal linkages, 
there are now famous bending flows (around diagonals) found by Kapovich and Millson 
\cite{KaMi-3D1996}, which are a special class of integrable Hamiltonian systems,
and the singularities of these systems have been described in a recent paper by
Bouloc \cite{Bouloc-3D2018}. We want to see to what extent these results can be 
generalized to other situations, especially planar linkages, and also 3D linkages 
which are more general than just polygons. It turns out that, for 2-dimensional 
(2D) polygons, we have a simple construction which may be viewed as an analogue 
of Kapovich-Millson systems: given a 2D $n$-gonal linkage ($n \geq 4$), just 
choose arbitrary $[(n-3)/2]$ ($[\cdot ]$ denotes the integer part of a real number) 
non-crossing diagonals which decomposes it into $[n/2] - 1$ quadrangles and 
at most one triangle (one if $n$ is odd and zero if $n$ is even). Then we get a 
non-Hamiltonian integrable system on the configuration space of type 
$([n/2] - 1,[(n-3)/2])$, i.e., with $[n/2] - 1$ commuting vector fields and 
$[(n-3)/2]$ common first integrals; ($[(n-3)/2] + [n/2] - 1 = n-3$ is the 
dimension of the configuration space). Each diagonal gives one first integral 
(its length, or rather its length squared to be smooth), each quadrangle has 
one degree of freedom and gives rise to one vector field, 
and they all commute with each other (see the rest of this paper for details). 
For example, when $n=7$, the configuration space of a generic planar heptagonal 
linkage is a 4-dimensional manifold, on which we have integrable systems of type 
(2,2), i.e., with two commuting vector fields and two common first integrals, 
by this construction. 

\begin{figure}[!ht]
\centering
\includegraphics[height = 0.36 \linewidth]{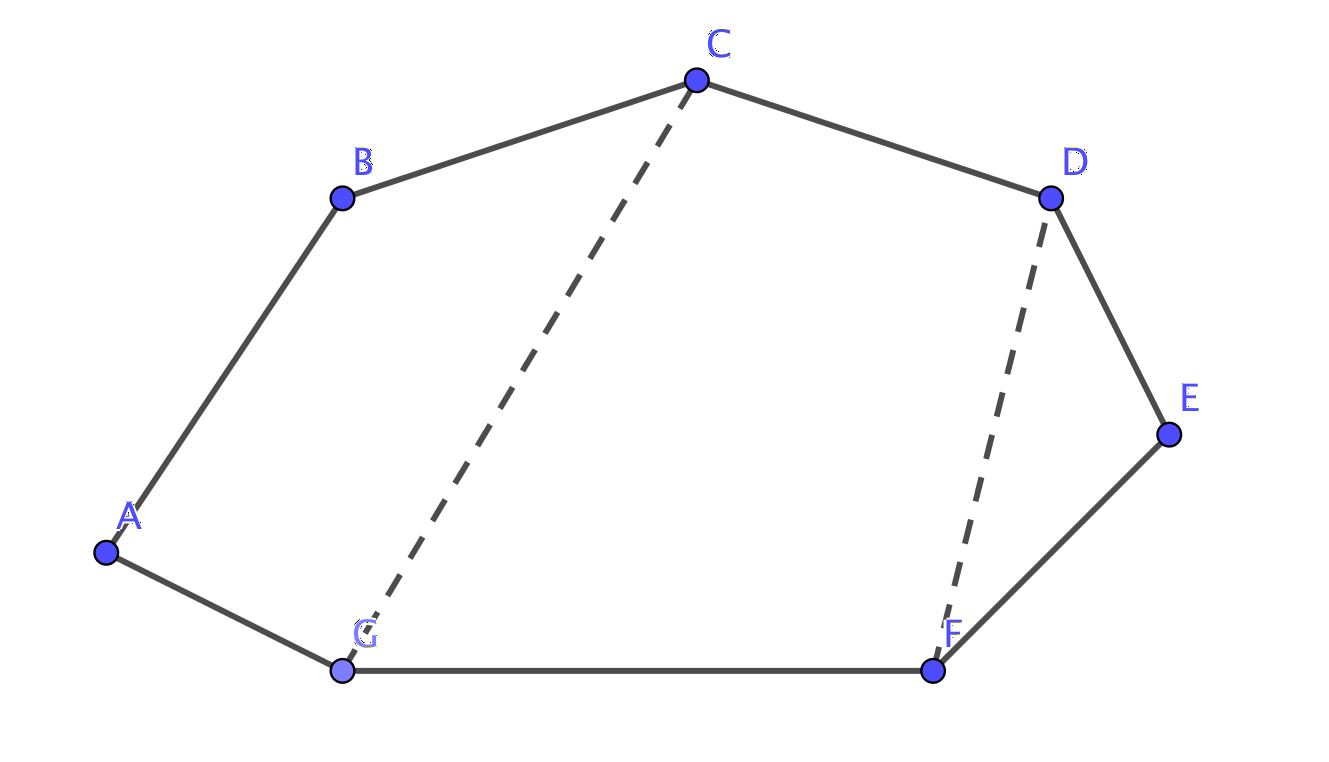} 
\caption{An integrable system on a the 4-dimensional configuration space  of a 
7-gon linkage: the two ``4-bar mechanisms'' $ABCG$ and $CGFD$ provide
two commuting vector fields  and the square lengths $|CG|^2$ and $|DF|^2$
are their smooth common first integrals.}
\label{fig:7gon}
\end{figure}

Recall that the configuration space of 3D linkages admits a natural symplectic 
structure \cite{Klyachko-3D1994,KaMi-3D1996,HaKn-3D1997}. However, for 2D linkages, 
a priori, we have neither a natural symplectic structure, nor Hamiltonian systems;
however, we have Nambu structures  and so our vector fields can be constructed 
in the spirit of Nambu mechanics \cite{Nambu-Nambu1973}.
Recall that  a Nambu structure 
of order $q$ on a manifold is a $q$-vector field which spans a (singular) integrable 
$q$-dimensional distribution and which may be viewed as a leaf-wise contravariant 
volume element on the corresponding singular $q$-dimensional foliation; see, e.g.,
\cite{DufourZung-PoissonBook2005,MinhZung-Foliations2013}.
Contravariant volume elements are particular cases of Nambu structures  and
they were used by Nambu in \cite{Nambu-Nambu1973} as a way to generalize Hamiltonian
mechanics, whence the name. 

$\bullet$ Our third motivation is given by the universality theorem of
Kapovich and Millson \cite{KaMi-Universal2002}, which was announced many years 
before by Thurston \cite{Thurston-Shapes1987} but without a full proof. This 
theorem states that, given any closed (compact without boundary) smooth manifold, 
there is a planar linkage whose configuration space is 
diffeomorphic to a disjoint union of a finite number of copies of that manifold. 
Thus, by studying natural integrable systems on configuration spaces of planar 
linkages, we can obtain interesting integrable non-Hamiltonian systems  on 
\textit{any} closed manifold. The theory of integrable non-Hamiltonian systems
(in the sense of Bogoyavlensky \cite{Bogo-Extended1998}) is relatively new 
compared to integrable Hamiltonian systems, especially when it comes to their 
topological and geometric properties (see \cite{Zung-NonHamiltonian2016} and 
references therein), so this huge class of integrable non-Hamiltonian systems 
coming from planar linkages will be very useful for the development of the theory. 
\medskip
We summarize below the main results of the paper.

The first result of this paper, presented in Section \ref{section:ConfigurationSpace},
is that every configuration space of planar linkages admits a natural volume form.
This volume form, together with a choice of diagonals in the linkage, gives rise
to Nambu structures on the configuration space. 

In Section \ref{section:Integrable} we give a simple general method for constructing
integrable systems on configuration spaces of planar linkages, by decomposing such 
spaces into cross  products of smaller configuration spaces, and lifting integrable 
systems from these smaller spaces. In particular, we introduce the notion of 
\textit{cross  products of integrable systems} in this section and, as a special 
case of this cross  product construction, we find, what we call, 
\textit{multi-Nambu integrable systems}.

In Section \ref{section:Final}, the last section of this paper, we list some open 
questions and final remarks. In particular, we mention that our first integrals 
are Bott-Morse functions on many configuration spaces, which are helpful in topology. 
For example, in the case of hexagons, by using these Bott-Morse
functions, one immediately gets the fact that the corresponding configuration
spaces belong to a special class of very well-studied 3-manifolds named 
graph-manifolds by Waldhausen \cite{Waldhausen-Graphs1967}; this
fact also fits well with Fomenko's topological theory of integrable systems
on 3-manifolds \cite{Fomenko-2DOF1989}. 

\section{Configuration spaces of planar linkages}
\label{section:ConfigurationSpace}

\subsection{Definition of the orbit and configuration spaces} 
\hfill

In this paper, a \textit{\textbf{planar linkage}} $\cL$   consists of:
\begin{itemize}
\item A finite graph $\cG = (\cV,\cE)$, where $\cV$ denotes 
the set of \textbf{\textit{vertices}}  and $\cE$ denotes 
the set of \textit{\textbf{edges}} of $\cG$.
For any two vertices $A,B \in \cV$ there is at most one edge connecting them,
denoted $AB$ or $BA$. An edge is also called a \textbf{\textit{bar}} 
in $\cL$.
\item Each edge has an associated positive length, i.e., there is a 
length function
$$ \ell : \cE \to \bbR_+:= \,]0, \infty[$$
which specifies the length $\ell(AB)$ of each edge $AB \in \cE$.
\item A subset $\cB \subset \cV$ of vertices, which may be empty, called 
the set of \textbf{\textit{fixed vertices}} or \textbf{\textit{based 
points}}, which comes together with a \textbf{\textit{fixing map}}
$$ \beta: \cB \to \bbR^2$$
that attaches $\cL$ to the Euclidean plane $\bbR^2$. The vertices in
$\mathcal{V} \setminus \mathcal{B} $ are called \textbf{\textit{movable}}.
\end{itemize}

The pair $ (\mathcal{G},\mathcal{B} )$, where the lengths and the fixing 
are not specified, is called the \textbf{\textit{linkage type}} of 
$ \mathcal{L} $.

A \textbf{\textit{realization}} of a planar linkage 
$\cL = (\cV,\cE, \ell,\cB,\beta)$ is a map
$$ \mathbf{r} : \cV \to \bbR^2$$
which respects the based points and the bar lengths. In other words, 
$$  |\mathbf{r} (A) \mathbf{r} (B)| = \ell(AB),$$
for any $A,B \in \cV$, and the restriction of $\mathbf{r} $ to the set 
$\cB$ of base vertices coincides with $\beta$. Of course, the length 
condition means that $\mathbf{r} $ can be extended to a map from the whole 
graph $\cG$ to $\bbR^2$, that maps each edge of $\mathcal{G}$  to a straight 
line segment of the same length. In a planar realization of a linkage, the 
bars can cross each other, i.e., we do not require the map from the 
graph to $\mathbb{R}^2$ to be injective. The crossing of bars 
may be realized in practice by laying them on parallel planes. 

The space of all realizations of a planar linkage $\cL$ is denoted by 
$\cO_\cL$; it is also called the \textbf{\textit{orbit space}} of $\cL$. 
Two different realizations $\mathbf{r} _1$, $\mathbf{r} _2$ of $\cL$ are 
said to be \textbf{\textit{isometric}} if there is an orientation and length 
preserving map $\phi: \bbR^2 \rightarrow  \mathbb{R}^2$ (i.e., 
$\phi \in SE(2)$, the Euclidean group) such that 
$\phi \circ \mathbf{r} _1 = \mathbf{r} _2$. The set of all 
realizations of $\cL$ modulo isometries is called the 
\textbf{\textit{configuration space}} and is denoted by $M_\cL$.

If $\cB = \emptyset$ is empty then the realizations of $\cL$ can move freely 
in $\bbR^2$; in this case $SE(2)$ acts on $\cO_\cL$, the action is free
if $ \mathcal{L} $ has at least one bar, and we have  
$$M_{\cL} = \cO_{\cL} / SE(2).$$
If $\cB$ consists of exactly one base point $B$, then the configurations 
cannot be translated freely in $\bbR^2$ but can turn around $B$; in this 
case we have
$$M_{\cL} = \cO_{\cL} / SO(2).$$
If, however, $\cB$ consists of at least two distinct base vertices, then two 
realizations of $\cL$ are isometric if and only if they coincide; in this 
case we have
$$M_{\cL} = \cO_{\cL}.$$
Let $G_ \mathcal{L} $ be the group $ \{e\}$, $SO(2)$, or $SE(2)$,
depending on whether $\mathcal{B}$ has at least two, one, or zero fixed vertices.
With this notation, in view of the considerations above, we have 
$M_\mathcal{L}=\mathcal{O}_\mathcal{L}/G_\mathcal{L}$.

Note that  when studying the configuration space of a linkage, we can 
choose between the version with at least 2 base points and the version 
without any base point. Indeed, if $\cL = (\cV,\cE, \ell,\cB,\beta)$, 
where $\cB$ has at least 2 points, then we can add all edges $BC$ with 
vertices $B,C \in \cB$, together with their lengths 
$\ell(BC) = |\beta(B)\beta(C)|$, then forget these base vertices, and we
get a linkage with empty base. This new linkage, now without any base 
vertices, has the same configuration space as the old linkage.
Conversely, if $\cL$ does not have any base points, then we can choose an
arbitrary bar in $\cL$ and fix it to $\bbR^2$, so that its vertices 
become base points, without changing the configuration space of $\cL$. 

\begin{figure}[!ht]
 \centering
\includegraphics[width = 0.7\linewidth]{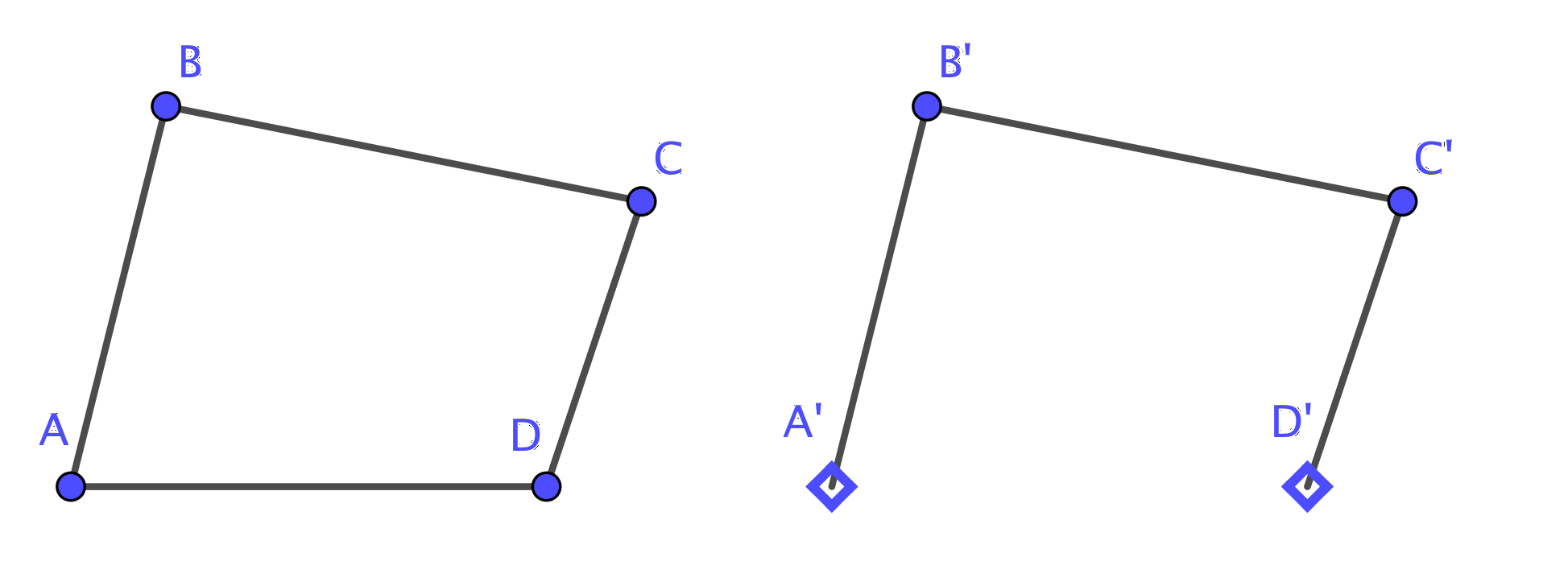} \\
\caption{Two versions of a quadrangle (with and without 
base points), also called a 4-bar mechanism, with the 
same configuration space.}
\label{fig_quadrangle}
\end{figure}

For example, consider a \textbf{\textit{quadrangle}} linkage $\cL$ (also 
called a \textbf{\textit{4-bar mechanism}}) consisting of four vertices 
$A,B,C,D$ and  four edges with four fixed lengths $\ell_1 = |AB|$,  
$\ell_2 = |BC|$,  $\ell_3 = |CD|,\ell_4 = |DA|$, and no base point. 
Consider the same quadrangle, but now with two base points $A$ and $D$, 
and only three free edges $AB$, $BC$, $CD$ (the edge $DA$ can be removed 
because $A$ and $D$ are already fixed). We get a new linkage $\cL_2$, 
whose configuration space is exactly the same as the configuration space 
of $\cL_1$ (see Figure \ref{fig_quadrangle}).

\begin{remark}{\rm
When talking about a linkage $\cL$, one often assumes that it is connected 
(i.e., the corresponding graph is connected). However, for our study, it 
is convenient to allow linkages which are disconnected. It is clear that 
if $\cL$ is a disjoint union of two linkages $\cL_1$ and $\cL_2$, then
its base is also the disjoint union of $\cB_1$ and $\cB_2$ and
its orbit space is simply the direct product of the two corresponding 
orbit spaces:
$$ \cO_{\cL_1 \sqcup \cL_2} = \cO_{\cL_1} \times \cO_{\cL_2}.$$

The formula for the configuration space of ${\cL_1 \sqcup \cL_2}$ is
more complicated. Namely, we have six  cases:

(i) If both $\mathcal{B}_1$ and $\mathcal{B}_2$ have at least two 
vertices, then $\mathcal{O}_\mathcal{L} = M_\mathcal{L} =
\cO_{\cL_1} \times \cO_{\cL_2} = M_{\mathcal{L}_1} \times  
M_{\mathcal{L}_2} $.

(ii) If $\mathcal{B}_1$ has at least two vertices and $\mathcal{B}_2$
consists of exactly one vertex, then $\mathcal{O}_\mathcal{L} = 
M_\mathcal{L} =\cO_{\cL_1} \times \cO_{\cL_2}$ but 
$M_{\mathcal{L}_1} \times M_{\mathcal{L}_2}= M_ \mathcal{L}/SO(2)$,
the $ SO(2)$-action being on the second factor. 

(iii) If $\mathcal{B}_1$ has at least two vertices and $\mathcal{B}_2
= \emptyset$,  then $\mathcal{O}_\mathcal{L} = M_\mathcal{L} =
\cO_{\cL_1} \times \cO_{\cL_2}$ but 
$M_{\mathcal{L}_1} \times M_{\mathcal{L}_2}= M_ \mathcal{L}/SE(2)$,
the $ SE(2)$-action being on the second factor.

(iv) If both $\mathcal{B}_1$ and $\mathcal{B}_2$ have each exactly one
vertex, then $\mathcal{O}_\mathcal{L} = M_\mathcal{L} =
\cO_{\cL_1} \times \cO_{\cL_2}$ but 
$M_{\mathcal{L}_1} \times M_{\mathcal{L}_2}= 
M_ \mathcal{L}/(SO(2)\times  SO(2))$, each $SO(2)$ acting separately
on its corresponding factor. 

(v) If $\mathcal{B}_1$ consists of one vertex and $\mathcal{B}_2= 
\emptyset$, then $M_\mathcal{L} = \mathcal{O}_\mathcal{L}/SE(2) =
\mathcal{O}_{\mathcal{L}_1} \times_{SO(2)}\mathcal{O}_{\mathcal{L}_2}$,
so we get fibrations $M_\mathcal{L} \rightarrow\mathcal{O}_{\mathcal{L}_1}$,
$M_\mathcal{L} \rightarrow\mathcal{O}_{\mathcal{L}_2}$ whose fibers are
$\mathcal{O}_{\mathcal{L}_2}$ and $\mathcal{O}_{\mathcal{L}_1}$, 
respectively. Thus we get an $SE(2)$-principal fiber bundle $M_\mathcal{L} 
\rightarrow M_{\mathcal{L}_1} \times M_{\mathcal{L}_2}$. The 
$SE(2)$-action is induced from the following $SE(2)$-action on the linkage:
fix a second vertex on $\mathcal{L}_1$ and let $SE(2)$ act on 
$\mathcal{L}_2$. The resulting action on $M_\mathcal{L}$ is free and the 
quotient $M_\mathcal{L}/SE(2)$ is $M_{\mathcal{L}_1} \times 
M_{\mathcal{L}_2}$.  

(vi) If $\mathcal{B}_1= \emptyset$ and $\mathcal{B}_2= \emptyset$,
we proceed as in the previous case, except that we now fix two vertices
in $ \mathcal{L} _1$. We get an $SE(2)$-principal fiber bundle 
$M_\mathcal{L} \rightarrow M_{\mathcal{L}_1} \times M_{\mathcal{L}_2}$. 
  
}
\end{remark}

\subsection{The dimension}
\hfill

Consider a linkage $\cL$ with $s$ fixed vertices $B_1,\ldots, 
B_s$ and $n$ movable vertices $A_1,\hdots, A_n$. If there are no bars, 
then each $A_i$ can take any position in $\bbR^2$ and we have
$$\cO_{(\mathcal{E},\mathcal{B},\beta)}\cong \bbR^2 \times \dots \times 
\bbR^2 \cong \bbR^{2n};$$
$\cO_{(\mathcal{E},\mathcal{B},\beta)}$ denotes the set of all realizations
of the planar linkage $\mathcal{L}$ but with variable lengths of the bars. 
If the linkage has $q$ bars, then we have a 
\textit{\textbf{joint length map}}
$${\boldsymbol{\mu}} = (\mu_1,\hdots,\mu_q): 
\cO_{(\mathcal{E},\mathcal{B},\beta)} \to \bbR_+^q,$$
where $\mu_i$ is half the length square of bar number $i$ in $\cL$.  
(We take the square to have smooth functions.) Clearly this map descends
to the quotient $ M_ {(\mathcal{E},\mathcal{B},\beta)} = 
\mathcal{O}_{(\mathcal{E},\mathcal{B},\beta)}/G_ \mathcal{L} $; 
we denote it by the same symbol ${\boldsymbol{\mu}}$, because there
is no danger of confusion.

For each fixed value of a $q$-tuple of lengths $(\ell_1,\hdots,\ell_q)$ 
we have
$$\mathcal{O}_\cL = {\boldsymbol{\mu}}^{-1}(\ell_1^2/2,\hdots,\ell_q^2/2) 
= \{L \in\mathcal{O}_{(\mathcal{E},\mathcal{B},\beta)} \mid 
\mu_i(L) = \ell_i^2/2, \ \forall 1  \leq i \leq q\}, $$
and a similar formula for $M_\mathcal{L}$.

A bar in a linkage type $(\mathcal{G},\mathcal{B})$  is called 
\textit{\textbf{redundant}} if its length is dependent on the lengths 
of the other bars. In other  words,  edge $i$ is redundant if the 
function $\mu_i$ is functionally dependent on the functions 
$\mu_1,\ldots,\mu_{i-1},\mu_{i+1},\ldots,\mu_q$. For example, 
consider a complete link of 4 vertices and 6 edges. Then one of the edges 
(it does not matter which one) is redundant. Eliminating a  
redundant edge does not change the dimension of the configuration space. 
So, to compute the dimension of  the configuration space of a linkage $\cL$, 
we may assume that its linkage type $(\mathcal{G},\mathcal{B})$ has all bars 
non-redundant. In this case, we have
$$ \dim M_\cL  = 2n - q$$
for any regular (generic) value of the $q$-tuple of edge lengths. 
(This formula is essentially the classical so-called Chebychev–Grübler–Kutzbach 
criterion for the degree of freedom of a kinematic chain; see, e.g.,
\cite{McCarthy-Designs2011}.)
In particular, 
 $$ q \leq 2n$$
in a linkage all of whose edges are non-redundant, 
i.e., the number of edges (bars) is at most twice the number 
of movable vertices. A linkage all of whose bars are non-redundant
is called a \textbf{\textit{non-redundant linkage}}. 

\subsection{The contravariant volume element}
\hfill

We want to define a natural volume multivector field on $M_\mathcal{L} $
(a ``contravariant volume", also called a Nambu structure of top order). 
The strategy is to construct it from 
a volume multivector field on $\mathcal{O}_\mathcal{L}$ and
on $ G_ \mathcal{L} $ (recall that $M_\cL = \cO_\cL/G_\cL$). 

First, recall that on any orientable manifold $ M$, any choice of
volume form $ \Omega_M$ uniquely defines a top degree nowhere
vanishing multi-vector field $ \Lambda _M$ by requiring it to
satisfy $\Lambda _M\intprod \Omega _M = 1$. We call such a 
top degree multi-vector field a \textbf{\textit{contravariant
volume}}.

Second, we note that on any Lie group $G$ there always is a 
left-invariant contravariant volume. Indeed, if 
$\{\xi_1, \ldots, \xi_p\}$ is a basis of its Lie algebra $\mathfrak{g}$, 
extend these vectors to left invariant vector fields $\xi_1^L, \ldots, 
\xi_p^L$ and form the contravariant volume 
$\xi_1^L\wedge \ldots\wedge \xi_p^L$.

Third, if the Lie group $G$ acts locally freely  on a dense set of
the manifold $M$ (i.e., all its isotropy subgroups at the points 
on this dense sent are discrete), there is
always a Nambu structure  of degree $ p=\dim G$ on $M$. Indeed,
if $\{\xi_1, \ldots, \xi_p\}$ is a basis of $\mathfrak{g}$, form
the corresponding infinitesimal generator vector fields $X_1,\ldots,X_p$;
then $\Pi_\mathcal{L} : = X_1\wedge \ldots \wedge X_p$ is a Nambu
structure on $M$ tangent to the generic $G$-orbits.
\medskip 

We return now to our problem of constructing a contravariant volume  
on $M_\mathcal{L}$. The Lie group $G_\mathcal{L}$ equals $\{e\}$, 
$SO(2)$, or $SE(2)$, depending on whether 
$ \mathcal{B} $ contains more than two, one, or zero vertices. Let 
$\Lambda_{G_\mathcal{L}}$ be the contravariant volume on  
$G_\mathcal{L}$ which equals hence\\ 
$\bullet$ $\Lambda_{G_\mathcal{L}}=1$ if $G_\mathcal{L}=\{e\}$,\\
$\bullet$ $\Lambda_{G_\mathcal{L}}=Z=a \frac{\partial }{\partial b} -
b \frac{\partial }{\partial a} $, the infinitesimal generator of
counterclockwise rotation in the plane $ \mathbb{R} ^2 =\{(a,b)\mid 
a, b \in  \mathbb{R} \}$, if $G_\mathcal{L}= SO(2)$,  \\
$\bullet$ $ \Lambda_{G_\mathcal{L}}=Z \wedge \frac{\partial }{\partial x} 
\wedge \frac{\partial }{\partial y}$, where $(x, y)$ are the coordinates of
$ \mathbb{R} ^2 $ in $ SE(2)=SO(2) \ltimes \mathbb{R} ^2 $.\\ 
This contravariant volume is given in terms of the obvious basis
elements of the Lie algebra of $G_\mathcal{L}$ and, in view of the 
previous discussion, since $G_\mathcal{L}$ acts freely on a dense subset 
of $\mathcal{O}_{\mathcal{L},\ell}$,  
produces a Nambu multivecor field $\Pi_\mathcal{L} $ on 
$\mathcal{O}_\mathcal{L}$, tangent to the generic 
$G_\mathcal{L}$-orbits on $\mathcal{O}_\mathcal{L}$. 

Next, we construct a contravariant volume on $\mathcal{O}_\mathcal{L} $.
Let $\mu _1, \ldots, \mu _q: \mathcal{O}_\mathcal{L} \rightarrow  
\,]0,\infty[$ be the collection of the squared length maps given by 
$ \mu _i( \mathbf{x} ):= \ell_i^2/2$, where $\ell_i$ is the length
of the edge $i$ and $ \mathbf{x} =((x_1, y_1), \ldots(x_n, y_n)) \in  
\mathcal{O}_\mathcal{L} $. Define
$$\Lambda _{\mathcal{O}_\mathcal{L}} := 
\left({\rm d} \mu _1 \wedge \ldots \wedge {\rm d} \mu_q \right) \intprod 
\left( \frac{\partial}{\partial x_1} \wedge
\frac{\partial}{\partial y_1}\wedge \cdots \wedge 
\frac{\partial}{\partial x_n} \wedge \frac{\partial}{\partial x_1}
\wedge \frac{\partial}{\partial y_n} \right). $$
This contravariant volume is unique up to the ordering of the bars 
in the linkage. Note that $\Lambda _{\mathcal{O}_\mathcal{L} }$ is 
invariant under the $G_ \mathcal{L} $-action.

Finally, let $ \Omega_{\mathcal{O}_\mathcal{L}}$ be the volume form on
$\mathcal{O}_\mathcal{L}$ dual to $\Lambda_{\cO_\cL}$,
i.e., $\Lambda _{\mathcal{O}_\mathcal{L} }\intprod 
\Omega_{\mathcal{O}_\mathcal{L}}=1$. Then $\Pi_\mathcal{L} \intprod 
\Omega_{\mathcal{O}_\mathcal{L}}$ is a basic form
because of $G_ \mathcal{L} $-invariance of the objects in this expression. 
So this form projects to a form $ \Omega_{M_\mathcal{L}}$ on  
$M_\mathcal{L}$. Let $\Lambda_{M_\mathcal{L}}$ be its 
dual multivector field, i.e., $\Lambda_{M_\mathcal{L}}\intprod 
\Omega_{M_\mathcal{L}} =1$. This is the volume multivector field 
on $M_\mathcal{L}$. We have proved the first part of the following  
statement.

\begin{theorem}
{\rm (i)} For each non-redundant linkage  $\cL$, its orbit 
 $\cO_\cL$ and configuration $M_\cL$ spaces admit 
natural canonical contravariant volumes $\Lambda_{\cO_\cL}$ 
and $\Lambda_{\cL }$, respectively.
These canonical contravariant volumes are uniquely determined by a 
choice of ordering of the edges of of $\cL$: if we change the ordering 
by a permutation $\sigma$ of parity $|\sigma|$ then these volume 
multivector fields are   multiplied by $(-1)^{|\sigma|}$.
 
{\rm (ii)} Consider two vertices $C,D$ of a non-redundant linkage $\cL$. 
Denote by $\mu_1 = |CD|^2/2$ the half squared $CD$-length function on 
$\cO_\cL$ and $M_\cL$. Assume that $\mu_1$ is not a constant function, 
so that the linkage $\cL_1$  obtained from $\cL$ by attaching to it 
an additional bar $CD$ with the given length $\ell(CD) = c$ is 
non-redundant. Then $\cO_{\cL_1}$ (resp., 
$M_{\cL_1}$) is the 
level set given by $\mu_1 = c^2/2$ in $\cO_{\cL}$ (resp., 
$M_{\cL}$) and we have $\dim \mathcal{O}_{\mathcal{L}_1} = 
\dim \mathcal{O}_ \mathcal{L} -1$, $\dim M_{\mathcal{L}_1} = 
\dim M_ \mathcal{L} -1$, and 
$$ \Lambda_{\cO_{\cL_1}} = {\rm d} \mu_1 \intprod \Lambda_{\cO_{\cL}} 
\quad \text{and}\quad \Lambda_{\cL_1} = {\rm d} \mu_1 \intprod \Lambda_{\cL}$$
restricted to these level sets.
\end{theorem}

The proof of (ii) is a direct verification.

\begin{figure}[!ht]
\centering
\includegraphics[width = 0.7\linewidth]{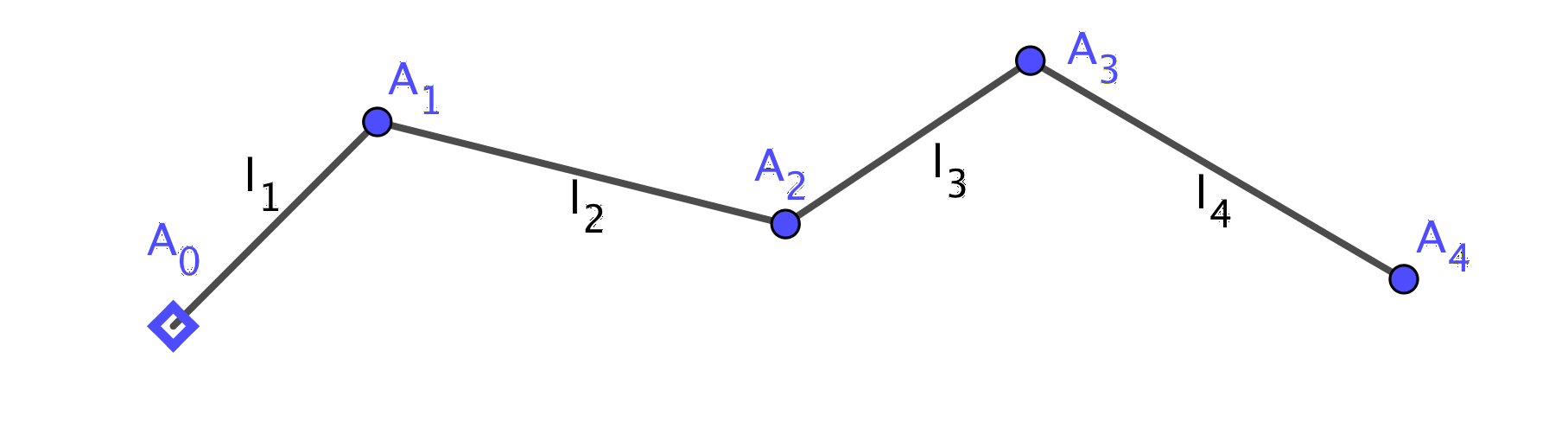} \\
\caption{The snake: the canonical contravariant volume elements
on its orbit space and configuration space are independent of its lengths.}
\label{fig:Snake}
\end{figure}

\begin{example} (The snake). \label{example:Snake}
{\rm
Consider the linkage $\cL$ which consists of
one base vertex $A_0$ fixed at the origin of $\bbR^2$ ($A_0 = (0,0)$), 
$k$ movable vertices $A_1 = (x_1,y_1) ,\hdots,A_k = (x_k,y_k) \in \bbR^2$, 
and $k$ fixed lengths $|A_{i-1}A_i| = c_i > 0$ ($i=1,\hdots, k$). 
Then it is easy to check that
$\cO_\cL \cong \mathbb{T}^{k}$ with periodic coordinates 
$\theta_i = \arctan((y_{i} -y_{i-1})/(x_{i} -x_{i-1})) \pmod {2\pi}$,
and the canonical contravariant volume element on $\cO_\cL$ is
$$ \Lambda_{\cO_\cL} = \dfrac{\partial}{\partial \theta_1}
\wedge \hdots \wedge \dfrac{\partial}{\partial \theta_k},$$
which does not depend on the lengths $c_1,\hdots,c_k$. The
configuration space $M_\cL$ is isomorphic to 
$\mathbb{T}^{k-1}$ with periodic coordinates
$\tilde\theta_i = \theta_{i+1} - \theta_i \pmod {2\pi}$, and its
canonical contravariant volume element is (up to a sign)
$$ \Lambda_{\cL} = \dfrac{\partial}{\partial \tilde\theta_1}
\wedge \hdots \wedge \dfrac{\partial}{\partial \tilde\theta_{k-1}}.$$
}
\end{example}

\begin{figure}[!ht]
\centering
\includegraphics[height = 0.38\linewidth]{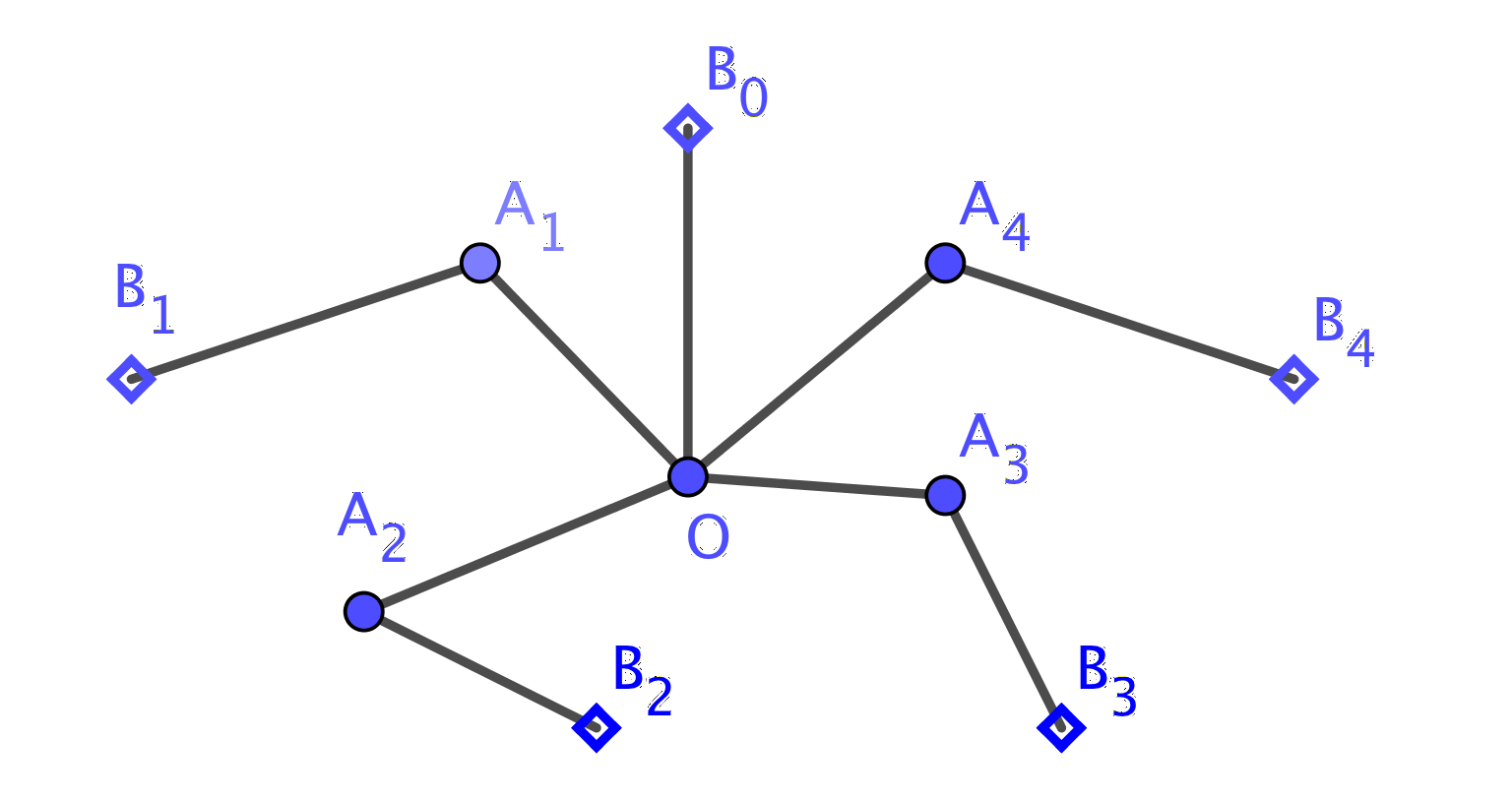} 
\includegraphics[height = 0.38 \linewidth]{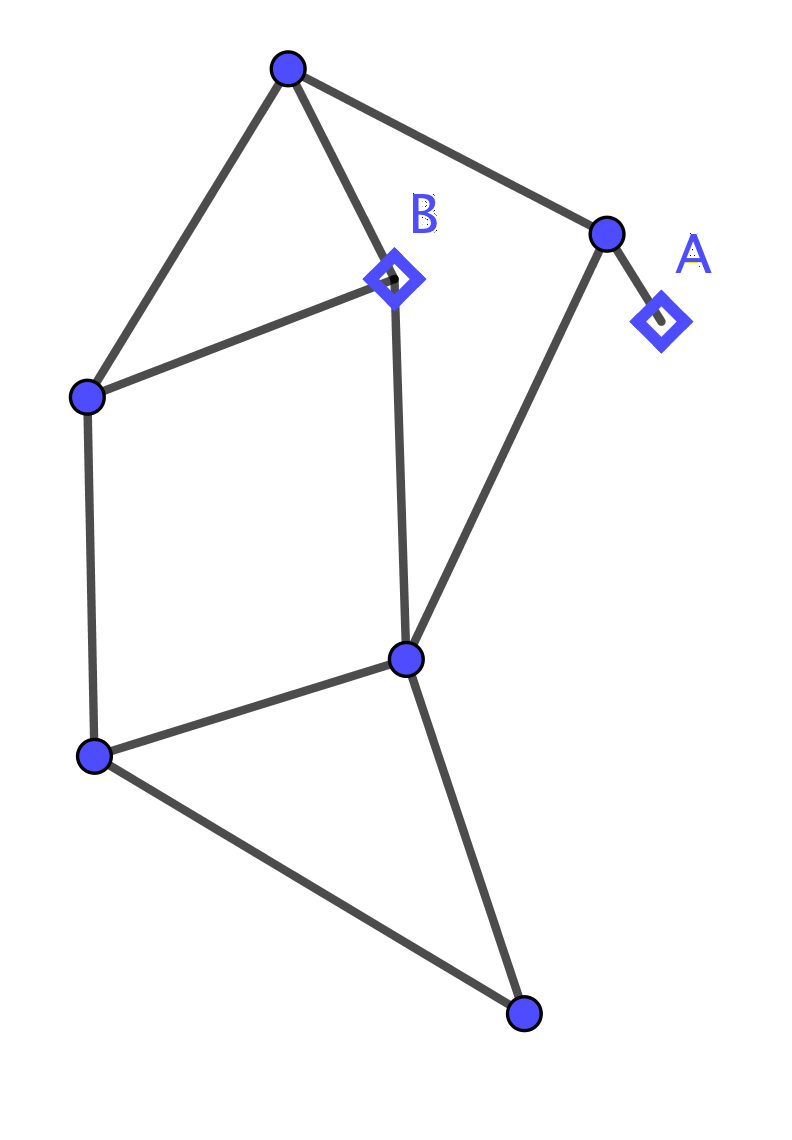} 
\\
\caption{Left: An 1DOF spider, with a fixed ``head'' $B_0$, 
fixed ``feet'' $B_1,\hdots,B_n$, ``body'' $O$, and ``legs'' 
$A_iB_i$ ($i=1,\hdots,n$). Right: ``Strandbeest'', invented 
by the artist Theodorus Jansen in 1990s, which can be used 
as legs in mechanical animals.}
\label{fig:Spider}
\end{figure}

\begin{example} (Curve-drawing linkages)
{\rm
Linkages with one degree of freedom (1DOF for short) 
are those whose configuration spaces have dimension 1:
generically, such configuration spaces are a finite union 
of closed curves. A planar 1DOF linkage is also called
\textit{\textbf{curve-drawing linkage}}  
because, after attaching it to a plane at the fixed vertices 
(if there are no fixed vertices then one can choose a bar 
and fix its two ends without changing the configuration 
space of the linkage), a movable vertex on the linkage will  
draw an algebraic curve when the linkage moves. By a famous 
en the linkage moves. By a famous 
1876 theorem of Kempe  \cite{Kempe-Curves1876}, any 
algebraic curve can be obtained this way by some 1DOF planar 
linkage. Simple modern proofs of this theorem can be found 
in many places; see, for example, 
\cite{Farber-Robotics2008,JoSt-Configuration1999}. 
Well-known examples of 1DOF linkages include  quadrangles, 
spiders, the Strandbeest (Figure \ref{fig:Spider}), and so on. 
A contravariant volume element on a manifold of dimension 1 is 
simply a vector field. Thus, our construction shows that  on 
the configuration space of each curve-drawing linkage there is a unique
(up to orientation) natural canonical vector field. Such vector fields
will be important in our construction of integrable systems on general
configuration spaces of planar linkages.
}
\end{example}

\begin{example} (Zero-dimensional space).
{\rm
Consider a linkage $\cL_{ABC}$ formed by a triangle $ABC$ with 3 fixed 
lengths (i.e., 3 bars) $|AB| = c$, $|BC| = a$, $|CA| = b$, without base 
points. This linkage is, of course, rigid, so its configuration space 
is just two points (corresponding to two possible orientations of the 
triangle on the plane). The contravariant volume element is just a 
number for each point on the configuration space. One of the numbers is
twice the area of the triangle $ABC$  and the other is its opposite 
(the same absolute value, but with negative sign). In particular, in 
the degenerate case, when the three points $A,B,C$ are aligned, 
the contravariant volume element vanishes. If we fix $A,B$ and 
erase the bar $AB$ from $\cL_{ABC}$, then it becomes another rigid 
linkage $\cL_{ABC}'$, with the same configuration space and the same 
contravariant volume element up to a natural isomorphism.
Now add to $\cL_{ABC}$ or $\cL_{ABC}'$ a vertex $D$ and a bar $CD$ 
with some length $|CD| = d > 0$. The configuration space of the new 
linkage is now two circles, and the corresponding contravariant 
volume element is $\pm 2S \dfrac{\partial}{\partial \theta}$, where
$\theta$ $\pmod {2\pi}$ is the angle formed by the vector 
$\overrightarrow{CD}$ with $\overrightarrow{CA}$, and $S$ is the 
area of the triangle $ABC$. 
}
\end{example}

The above simple example illustrates the following two phenomena:

$\bullet$ (\textit{Effect of fixing a base}) If $\cL'$ is obtained from 
$\cL$ by fixing two vertices $A,B$ of a bar $AB$ and then removing that 
bar, where $\cL$ is a connected linkage without base points, then not 
only are $M_\cL$ and $M_{\cL'}$ naturally isomorphic, but under this 
isomorphism we also have $\Lambda_\cL = \Lambda_{\cL'}$. (This fact is 
not surprising, because fixing $A,B$ amounts to contracting with the 
volume element of $SE(2)$ in the formula for $\Lambda_{\cL}$).

$\bullet$ (\textit{Effect of homothety}) If we multiply the length of 
every edge of a linkage $\cL$ by a positive constant $c$, then we get 
a linkage $c\cL$ homothetic to $\cL$, and $M_{c\cL}$ is naturally 
isomorphic to $M_{\cL}$. However, under this natural isomorphism we 
have the following homothety formula, whose proof is a straightforward 
verification:
\begin{equation}
\Lambda_{c\cL} = c^{2(|\cE| - (|\cV| + |\cB|))} \Lambda_{\cL}
\end{equation}
where $|\cV| - |\cB|$ is the number of movable vertices and $|\cE|$ is 
the number of edges, in the case when  $|\cB| \geq 2$ (i.e., when 
$M_\cL = \cO_\cL$) and $\cL$ is connected. In the case when 
$\cB = \emptyset$ and $\cL$ is connected, the above formula must 
be modified as follows: 
\begin{equation}
\Lambda_{c\cL} = c^{2(|\cE| - |\cV| + 1)} \Lambda_{\cL}
\end{equation}

\subsection{Nambu structures for linkages with marked diagonals} \hfill 

Consider a linkage $\cL$. If two vertices $A,B$ of $\cL$
are not connected by a bar of $\cL$, then we say that
$AB$ is a \textbf{\textit{diagonal}} in $\cL$. 
Choose a set $\cD$ of diagonals of $\cL$
(which may be empty) and mark the elements of $\cD$; in
this case we say that they are 
\textbf{\textit{marked diagonals}}.

Given a linkage $\cL$ with marked diagonals
$\cD = \{A_1B_1,\hdots,A_kB_k\}$ ($k \geq 0$), we can 
define the following Nambu structure on the configuration 
space $M_\cL$:
\begin{equation}
\Lambda_{\cL,\cD} = ({\rm d} \mu_1 \wedge \hdots \wedge
{\rm d} \mu_k) \intprod \Lambda_\cL,
\end{equation}
where $\mu_i = |A_iB_i|^2/2$. (When $\cD = \emptyset$ is 
empty then, of course, $\Lambda_{\cL,\emptyset} = \Lambda_\cL$). 
We will use such Nambu structures for the construction of our 
integrable systems.

\section{Integrable systems on linkage spaces}
\label{section:Integrable}

\subsection{Cross  products of integrable systems} \hfill

In this subsection, we present a general method for constructing 
new integrable systems from other integrable systems. We call
this technique the \textbf{\textit{cross  product}} construction.

We begin by recalling the notion of non-Hamiltonian integrable
systems (\cite{Bogoyavlenskij-Integrable1996, Zung-NonHamiltonian2016};
for an Action-Angle Theorem in this setting, see \cite{Zung-AA2018}).
A vector field $X$ on an $n$-dimensional manifold $M$ is called 
\textbf{\textit{integrable}} if there exist $p$ vector fields 
$X_1=X$, $X_2,\ldots,X_p$ and $q$ functions $F_1,\ldots,F_q$ on 
$M$, $p+q=n$, such that $[X_i, X_j]=0$ for all $i,j=1, \ldots, p$,
$X_1\wedge\cdots\wedge X_p\neq0$ almost everywhere, 
${\rm d}F_\ell(X_i)=0$ for all $i =1, \ldots, p$ and $\ell=1, \ldots, q$,
and $\mathrm dF_1\wedge\cdots\wedge\mathrm dF_q\neq0$ almost everywhere.
Such an $n$-tuple $(X_1,\ldots,X_p,F_1,\ldots,F_q)$ is  
called an \textit{\textbf{integrable system of type $(p,q)$}}.
\medskip

Consider a family of smooth manifolds $M_i$ of dimensions $n_i$, 
where $i$ belongs to some finite index set $I$. We assume that 
there is another finite index set $K$ disjoint from $I$, such 
that for each $k \in K$ there is a nonempty subset $S_k 
\subset I$, and for each $i \in S_k$ there is a given smooth function
$f_{i,k}: M_i \to \bbR$ on $M_i$. Moreover, we also assume that
for each $i \in I$ there is an integrable non-Hamiltonian system
whose functionally independent first integrals are  
$f_{i,k}$ for $k \in K_i$, where $K_i = \{k\in K\ |\ i \in S_k\}$, 
and whose linearly independent commuting vector fields 
are $X_{i,1},\hdots,X_{i,p_i}$. (By definition, we have 
$p_i + |K_i| = n_i$). 

Consider now the \textbf{\textit{diagonal manifold}}
\begin{equation}
M = \Big\{ {\bf x} =(x_i) \in \prod_{i \in I} M_i \ \Big | \  
f_{ik}(x_i) = f_{jk} (x_j)\ \forall i,j\in I, \  \forall 
k \in S_i \cap S_j \subset  K\Big\}.
\end{equation}

Denote by $\pi_j: \prod_{i \in I} M_i \to M_j$ the projection maps. By 
abuse of notation, we denote the  restriction of $\pi_j$ to $M$ again 
by $\pi_j$, $j\in  I$. We assume that the functions 
$\pi_i^*f_{ik}- \pi_j^*f_{jk}$ ($k \in S_i \cap S_j$) have linearly 
independent differentials on $M$, except at the points satisfying
the obvious relations of the type  $\pi_i^*f_{ik} - \pi_j^*f_{jk} = 
(\pi_i^*f_{ik} - \pi_\ell^*f_{\ell k}) + (\pi_\ell^*f_{\ell k} - 
\pi_j^*f_{jk})$; therefore, $M$ is a smooth manifold of dimension 
$\dim M = \sum_{i \in I} n_i - \sum_{k \in  K} (|S_k| - 1) = 
\sum_{i \in I} n_i - \sum_{k\in  K} |S_k| + |K|$ (because each 
index $k\in  K$ creates $|S_k| - 1$ independent equations). Notice 
that $\sum_{k\in  K} |S_k|$ is the total number of functions 
$f_{ik}$  and that this number equals $\sum_{i \in I} (n_i-p_i)$, 
because for each $i$ we have exactly $|K_i|=n_i-p_i$ functions. 
Therefore,
\begin{equation}
\label{dim_M} 
\dim M = |K| + \sum_{i \in I} p_i.
\end{equation}
Denote by $\widehat{X}_{i,j}$ the horizontal lift of $X_{i,j}$ to 
$\prod_{i\in  I} M_i$ (for each $i\in I$ and each $j \leq p_i$). 
In other words, $\widehat{X}_{i,j}$ is the unique vector field on 
$\prod_i M_i$  whose projections by $\pi_j$ ($j\neq i$) are trivial 
and whose projection by $\pi_i$ is $X_{i,j}$. Obviously, the vector 
fields $\widehat{X}_{i,j}$ commute on $\prod_{i \in I} M_i$. 
Moreover, they admit the pull-backs of the functions $f_{ik}$ as 
common first integrals. Hence, they also admit the functions 
$\pi_i^*f_{ik} - \pi_j^*f_{jk}$ as first integrals  and are, 
therefore, tangent to $M$. 

For each $k\in K$, define $F_k: M \to \mathbb{R}$ by
\begin{equation}
\label{def_F_k} 
F_k({\bf x}) = f_{ik} (\pi(x_i))\quad \text{for some} \quad
i \in S_k.
\end{equation}
It is immediate from the definition of 
$M$ that $F_k$ is well-defined. Moreover, these functions 
$F_k$ are common first integrals of the commuting vector fields
$X_{i,j}$ ($i\in I, j \leq p_i$). Notice that the number of 
vector fields $X_{i,j}$ is $\sum_{i \in I} p_i$, the number 
of functions $F_k$ is $|K|$, and $|K|+\sum_{i\in  I} p_i$ is exactly 
the dimension of $M$ (see \eqref{dim_M}). Thus, we immediately 
obtain the following result.

\begin{theorem}
\label{thm:CrossedProductIntegrable}
With the above notations and assumptions, on the diagonal manifold 
$M$ of dimension $|K| + \sum_{i \in I} p_i$ we have an integrable 
system of type $\left( \sum_{i \in I} p_i, |K|\right) $, consisting of 
$\sum_{i \in I} p_i$ commuting  vector fields  $\widehat{X}_{i,j}$ 
($i \in I, j \leq p_i$) and  $|K|$ first integrals $F_k$ ($k \in K$)
given by \eqref{def_F_k}.
\end{theorem}

The above constructed integrable system on $M$ is called the
\textbf{\textit{cross product}} of the corresponding integrable systems
on $M_i, i \in I$.

\subsection{Multi-Nambu integrable systems} \hfill

Consider now the following special case of the general cross product 
construction in the previous subsection.

Assume that for each $i\in  I$, the set $\{f_{i,k} \mid k \in  K_i\}$ 
($K_i = \{ k \in K\mid i \in S_k\}$), of functionally independent 
functions  on $M_i$ has exactly $n_i-1$ elements,  where 
$n_i = \dim M_i$  and, moreover, each $M_i$ admits a Nambu structure 
$\Lambda_i$ of top order, i.e., a  $n_i$-vector field on $M_i$.
(In practice, $\Lambda_i$ is often the dual of a volume form on 
$M_i$, though the construction still works if $\Lambda_i$ has 
singular points.) Then, for each $i\in I$,  define the  
\textbf{\textit{Nambu vector field}} $X_i$ on $M_i$ by
the following contraction formula:
\begin{equation}
\label{eqn:Nambu}
X_i = \left( \wedge_{k\in  K_i} {\rm d} f_{i, k}\right)  
\intprod \Lambda_i. 
\end{equation}
(Choose an ordering on $K_i=\{k \in  K\mid i \in S_k \}$ to fix 
the sign of $\wedge_{k\in  K_i}{\rm d} f_{i, k}$).
Clearly, $X_i$ admits $n_i-1$ first integrals $f_{i,k}$ 
($k \in  K_i$), since $|K_i|=n_i-1$, 
and hence $X_i$ is a completely 
integrable vector field on $M_i$. In other words, on each $M_i$ 
we have an integrable system $(X_i, f_{i,k})$ of type 
$(1, n_i-1)$. So, in the notation of the previous subsection,
we have $p_i = 1$ and $\sum_{i\in  I} p_i = |I|$. Thus we get:

\begin{cor}
\label{cor:MultiNambu}
With the above notations and assumptions of this subsection, 
on the diagonal manifold $M$ of dimension
$|I|+|K| $ we have an integrable system of type $(|I|, |K|)$, 
consisting of $|I|$ commuting vector fields  
$\widehat{X}_{i}$ ($i \in I$) and $|K|$ first integrals 
$F_k$ ($k \in K$).
\end{cor}

In view of this construction, we  call the systems obtained 
in the corollary above \textbf{\textit{multi-Nambu integrable systems}}.

\begin{remark}(Separation of variables and symplectization). 

{\rm
i) The construction above not only gives integrable systems, but 
also a kind of separation of variables for such systems: the separated 
variables are on the different spaces $M_i$ ($i\in I$).  In classical 
mechanics, a Hamiltonian system admitting a complete separation
of variables, in the sense of Hamilton-Jacobi theory, is automatically 
integrable. Our systems are a special case of this general idea, even 
though, a priori, there is no natural (pre)symplectic or Poisson structure 
which turns them into Hamiltonian systems.

ii) An exception is the case when $\dim M_i = 2$, for every $i \in I$, and 
the Nambu structure $\Lambda_i$ on $M_i$ is non-degenerate, i.e., it is 
dual to a  volume form $\omega_i$ on $M_i$. In dimension 2, a volume form 
is the same as a symplectic form. In this situation we can pull back every 
$\omega_i$ and take their sum to form a \textit{presymplectic} form
$$ 
\omega = \sum_{i\in  I} \pi_i^* \omega_i 
$$
on the diagonal manifold $M$. Thus we obtain a presymplectic manifold 
$(M,\omega)$ together with an integrable Hamiltonian system 
$(\widehat{X}_i \ (i \in I),F_k\ (k \in K))$. Each vector field 
$\widehat{X}_i $ is Hamiltonian: $\widehat{X}_i  \intprod \omega = 
- {\rm d} F_k$, where $k \in K$ is the only
index in $S_i$ (i.e, there is only one function $f_{ik}$ on $M_i$). 
The rank of $\omega$ at a generic point on $M$ is $2|K|$ in this case.

iii) If $\dim M_i \geq 3$, then even if  $|I| = 1$, there is no  
natural way to turn the Nambu vector field $X_i$ on $M_i$ into 
a Hamiltonian vector field (with respect to a presymplectic,  
Poisson, or Dirac structure). The problem is that, in Nambu mechanics, 
we have $n-1$ functions which play equally important roles in 
the definition of the Nambu vector field, while in Hamiltonian 
mechanics we have just one function (the energy function). In 
order to turn a Nambu vector field into a Hamiltonian vector field, 
one needs to treat these functions differently and declare that 
among the $n-1$ function there is a special one which is more 
important than the other functions. Of course, this can be done. 
For example, a Nambu vector field 
$X = ({\rm d} F_1 \wedge \hdots \wedge {\rm d} F_{n-1}) \intprod \Lambda$, 
where $\Lambda$ is a Nambu structure of order $n$ on a manifold, 
may be viewed as a Hamiltonian vector field given by $F_1$ 
with respect to the Poisson structure 
$\Pi = ({\rm d} F_2 \wedge \hdots \wedge {\rm d} F_{n-1}) \intprod \Lambda$
of rank 2. Nevertheless, this is not very natural, especially 
in view of the diagonal construction.

iv) Moreover, there is no natural way to lift multi-vector 
fields from $M_i$ to the diagonal manifold $M$ (the horizontal 
lifting that we used for our vector fields does not work here 
because the result is not tangent to $M$). Thus, even if the 
$M_i$ are equipped with Poisson structures, $M$ is still not Poisson.
We cannot even lift Nambu structures $\Lambda_i$ from $M_i$ to $M$. 
So the Nambu structures in a multi-Nambu integrable system are 
not Nambu structures on the ambient manifold,  only on its 
projected components.
}
\end{remark}

\subsection{Decomposition of linkages and integrable systems}\hfill

We apply now the constructions in the previous subsections to
the configuration spaces of planar linkages.

\begin{figure}[!ht]
\centering
\includegraphics[width = 0.65\linewidth]{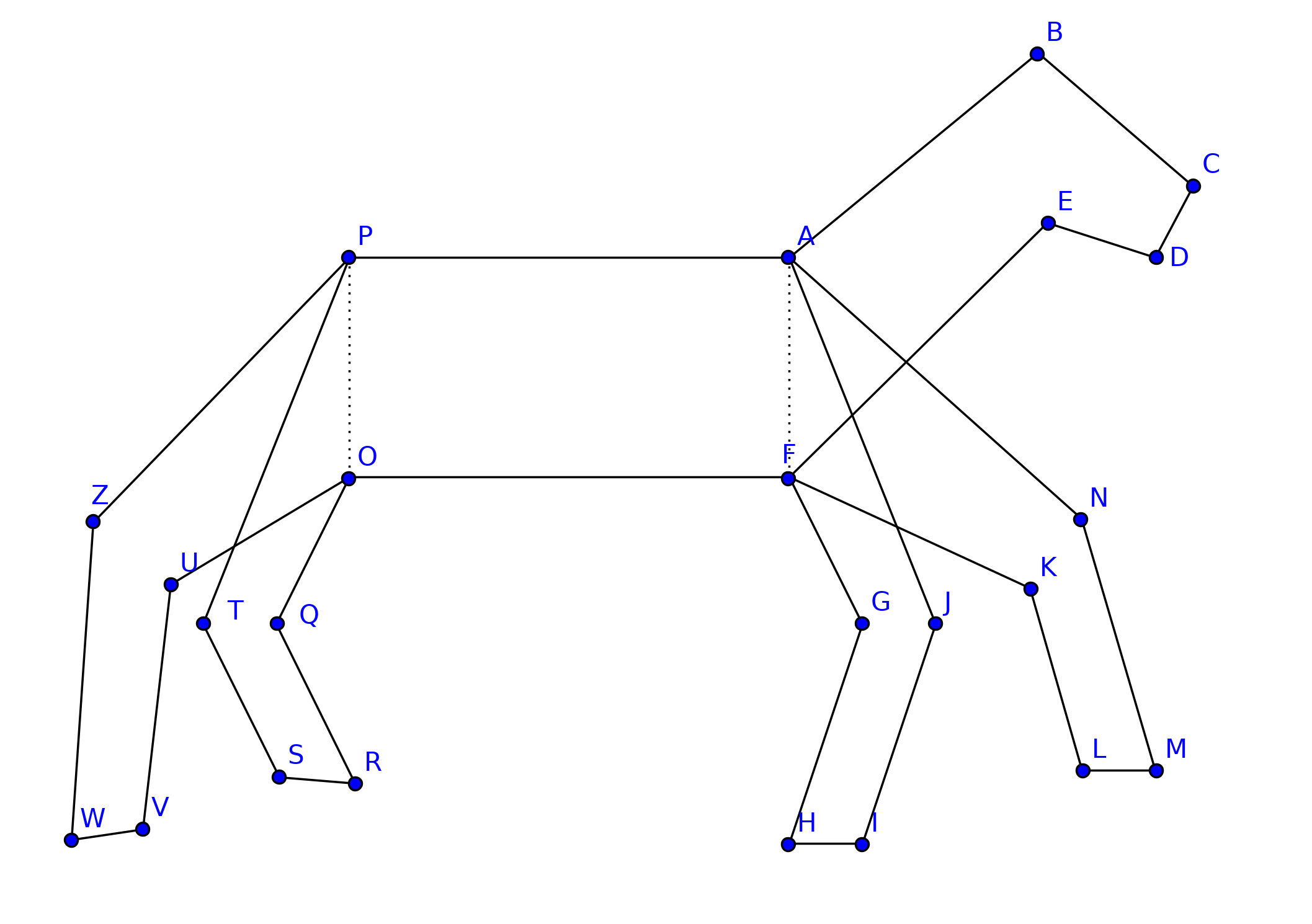} \\
\caption{Example of a connected sum with two joints 
and two marked diagonals $AF$ and $OP$: four 
sublinkages share the same joint $\{A,F\}$ and 
three sublinkages share the same joint $\{O,P\}$.}
\label{fig:horse}
\end{figure}

Consider a planar linkage $\cL = (\cV,\cE, \ell)$ without base  
points. (Recall that any non-trivial linkage with base points is 
equivalent to some linkage without base  points, configuration-wise). 
Choose a set $\cD$ of diagonals of $\cL$ such that $\cD$ decomposes 
$\cL$ into an \textit{acyclic semi-rigid connected sum} of a family 
of non-empty sublinkages $\cL_i = (\cV_i,\cE_i, \ell_i)$ ($i \in I$) 
with their corresponding sets $\cD_i$ of marked diagonals. This means 
that the following conditions are satisfied:

\begin{itemize}
\item For each $i\in I$, denote by $\cV_i$ the set of vertices of $\cL_i$. 
Then $\cL_i = (\cV_i,\cE_i, \ell_i)$ is the sublinkage of $\cL$ generated 
by $\cV_i$. This means that $\cE_i$ is the subset of $\cE$ consisting of
all edges of $\cL$ such that all of their vertices belong to $\cV_i$;
likewise, $\ell_i$ and $\cD_i$ 
are also defined in a natural way. (The diagonals in the family $\cD_i$
are the elements of $\cD$ whose vertices belong to $\cV_i$.) 

\item $\cV = \cup_{i\in I} \cV_i$, $\cE = \cup_{i\in I} \cE_i$, and
$\cD = \cup_{i\in I} \cD_i$.

\item There is a family $\cJ_k$ ($k \in K$) of sublinkages of $\cL$,
which are called \textbf{\textit{joints}} for the sublinkages $\cL_i$, 
in the following sense:  if $\cL_i \cap \cL_j \neq \emptyset$ for 
$i,j \in I$, then there is a $k \in  K$ such that 
$\cL_i \cap \cL_j = \cJ_k$. (It may happen that more than two 
sublinkges in the family $\cL_i$ ($i \in I$) share the same joint;
see Figure \ref{fig:horse}).

\item Every marked diagonal is a diagonal of a joint:
for any $VW \in \cD$ there is a joint $\cJ_k$ such that
$V,W$ are vertices of $\cJ_k$. 

\item The joints are disjoint: if $k_1 \neq k_2$ then $\cJ_{k_1}$ and
$\cJ_{k_2}$ don't have any common vertex.

\item \textit{Semi-rigidity}: For each $k \in K$, 
$\cJ_k$ has at least $2$ vertices. If we fix the lengths of the marked 
diagonals (from the family $\cD$) whose vertices lie in $\cJ_k$ 
then $\cJ_k$ becomes rigid, i.e., its configuration space becomes 
just one point (if not empty) once these diagonal lengths are fixed.

\item \textit{Acyclicity}: Consider the graph whose vertices are 
the elements of the index sets $I$ and $K$ (of course, these two 
index sets are assumed to be disjoint)  and whose edges are the 
couples $(i,k) \in I \times K$ (connecting the vertex $i$ to the 
vertex $k$) such that the sublinkage $\cL_i$ contains the joint 
$\cJ_k$. Then this graph is a tree (i.e., it has no cycles).
\end{itemize}

\begin{definition}
With the above notations and under the above conditions, we say 
that $(\cL,\cD)$ is the \textbf{acyclic semi-rigid connected sum} of
its sublinkages with marked diagonals $(\cL_i,\cD_i)$ ($i\in I$).
\end{definition}

\begin{theorem}
\label{thm:DecompositionIntegrableLinkages}
{\rm (i)} With the above notations, assume that $(\cL,\cD)$ is the 
acyclic semi-rigid connected sum of the family $(\cL_i,\cD_i)$ 
($i\in I$). Assume that the diagonal square length functions 
$\mu_{VW} = \ell(VW)^2/2$ ($VW \in \cD$)  are functionally 
independent  almost everywhere
on the configuration space $M_{\cL}$,
and that for every $i\in I$ the number of elements in $\cD_i$ 
is strictly less than the dimension $\dim M_{\cL_i}$ of the 
configuration space of $\cL_i$. Then $M_{\cL}$ is the cross 
product of the family $M_{\cL_i}$ with respect to the marked 
square length functions:
$$
M_\cL \cong \Big\{\mathbf{x}=(x_i) \in \prod_{i \in I} M_{\cL_i}\ | 
\ \mu_d(x_i) = \mu_d(x_j)\ \forall i,j\in I, \ \forall 
d \in \cD_i \cap \cD_j \Big\}.
$$
In particular, we have
$$ 
\dim M_\cL = \sum_{i \in  I} (\dim M_{\cL_i} - |\cD_i|) + |\cD|.
$$

{\rm (ii)} Assume, moreover, that on each configuration space 
$M_{\cL_i}$ we have an integrable system of type $(p_i,q_i)$, 
where $q_i = |\cD_i|$ and $p_i = \dim M_{\cL_i} - q_i$, given 
by $p_i$ commuting vector fields $X_{i,1},\hdots,X_{i,p_i}$ 
and $|\cD_i|$ diagonal squared length functions $\mu_d$ 
($d \in \cD_i$) which serve as first integrals for the system. 
Then by the method of cross product construction we get an 
integrable system of type $\left(\sum_{i\in  I} p_i, |\cD|\right)$ 
on $M_\cL$ given by the vector fields $\widehat{X}_{i,j}$ (which 
are horizontal liftings of the vector fields $X_{i,j}$, 
$i \in I$, $j \leq p_i$) and diagonal squared length functions 
$\mu_d$ ($d\in \cD$).

{\rm(iii)} In particular, if $\dim M_{\cL_i} - |\cD_i| = 1$ 
for every $i \in I$, then on each $\dim M_{\cL_i}$ we have 
an integrable system of type $(1,|\cD_i|)$ given by 
$|\cD_i|$ marked diagonal square length functions and the 
Nambu vector field $X_i$ given by formula \eqref{eqn:Nambu} 
with respect to these functions and the natural 
contravariant volume element on $\dim M_{\cL_i}$. Consequently, 
in this case on $M_\cL$ we have a multi-Nambu integrable 
system of type $(|I|,|\cD|)$.
 \end{theorem}

The proof of the above theorem is straightforward. In assertion 
(i), the semi-rigidity and the acyclicity conditions ensure that,
for each point $(x_i) \in \prod_i M_{\cL_i}$ satisfying 
the ``equal common lengths'' constraints, we can glue the 
configurations $x_i$ together along the joints to create 
a configuration $x \in M_{\cL}$, in a unique way up to isometries. 
Assertions (ii) and (iii) follow immediately from of assertion (i),
Theorem \ref{thm:CrossedProductIntegrable}, and Corollary
\ref{cor:MultiNambu}.

\begin{example}{\rm The case of $n$-gons.
Denote by $(A_1,...,A_n)$ the consecutive vertices of an $n$-gon.
The configuration space $M_{\cL}$ in this case has
$2(n-2) - (n-1) = n-3$ dimensions.

If $n$ is even, we can divide the $n$-gon $(A_1 \dots A_n)$ into
$(n/2) - 1$ quadrangles by adding $(n/2) - 2$ diagonals to the 
linkage; for example, $A_1A_3, A_1A_5, \hdots, A_1A_{2n-1}$.
Each quadrangle may be viewed as a curve-drawing linkage and 
thus produces one vector field. So we have an integrable system 
of type $\big((n/2) - 1, (n/2)-2\big)$, i.e., with $(n/2)-1$ 
commuting vector fields and $(n/2) - 2$ common first integrals.

If $n$ is odd,  divide the $n$-gon $(A_1 \dots A_n)$ into
$(n-3)/2 - 1$ quadrangles plus one triangle by adding 
$(n-3)/2$ diagonals to the linkage; 
for example, $A_1A_3, A_1A_5, \hdots, A_1A_{2n-1}$.
So we get  an integrable system of type $((n-3)/2, (n-3)/2)$.
}
\end{example}

\begin{remark}
{\rm
i) In assertion (iii) of Theorem \ref{thm:DecompositionIntegrableLinkages}, 
once the lengths of the marked diagonals are fixed, the sublinkages $\cL_i$ 
become curve-drawing linkages (i.e, linkages whose configuration spaces are
1-dimensional). This shows that curve-drawing linkages play an important 
role in the construction of integrable systems on the configuration spaces 
of more general planar linkages.

ii) In a sense, the above construction is quite similar to the construction 
of bending flows by Kapovich and Millson \cite{KaMi-3D1996} for obtaining 
integrable systems on configuration spaces of 3D polygons. 

iii) The case when then dimension of $M_{\cL_i}$ is equal to (instead of 
being greater than) the number of marked diagonals in ${\cL_i}$ for some 
$i \in I$ is excluded in Theorem \ref{thm:DecompositionIntegrableLinkages}
for simplicity. However, this case can, in fact, be treated similarly. In order 
to address this case, we have to work on the level of orbit spaces instead 
of configuration spaces, so it is a bit more complicated. Some components 
${\cL_i}$ for which the number of marked diagonals is equal to 
$\dim \cL_i $ still give rise to a non-trivial (and periodic) vector field 
in the final integrable system on $M_{\cL}$. For example, see the case of 
the snake (Example \ref{example:Snake}), which is cut into pieces.  Each 
piece has 0-dimensional configuration space, but the configuration space 
of the snake is a torus with natural rotational commuting vector fields. 

The case when a joint $J_k$ consists of just one point is also excluded
in Theorem \ref{thm:DecompositionIntegrableLinkages} for simplicity. 
However, this case can also be treated similarly. Such a $J_k$
increases the difference $M_\cL - \sum_i \dim M_{\cL_i}$ and creates 
additional rotational vector fields on $M_\cL$. 

Likewise, the case with base points, also excluded in Theorem
\ref{thm:DecompositionIntegrableLinkages}, can also be treated similarly; 
just put all the base points in one component. 

iv) There is an obvious method to create integrable systems on
configuration spaces, which works for any planar linkage.
If the dimension of the configuration space $M_\cL$ is $n$, 
then just mark $n-1$ diagonals whose length functions are independent 
on $M_\cL$. Then, together with the contravariant volume element on
$M_\cL$, we   immediately get a Nambu integrable system of
type $(1,n-1)$ whose first integrals are squared length functions
of the marked diagonals  and whose vector field is given by
formula \eqref{eqn:Nambu}. However, such ``obvious'' systems 
have just one vector field each, so if we want to obtain systems 
with more vector fields using the cross product method, we have to 
decompose the linkage $\cL$ into many components $\cL_i$, the more, 
the better. 

\begin{figure}[!ht]
\centering
\includegraphics[width = 0.6\linewidth]{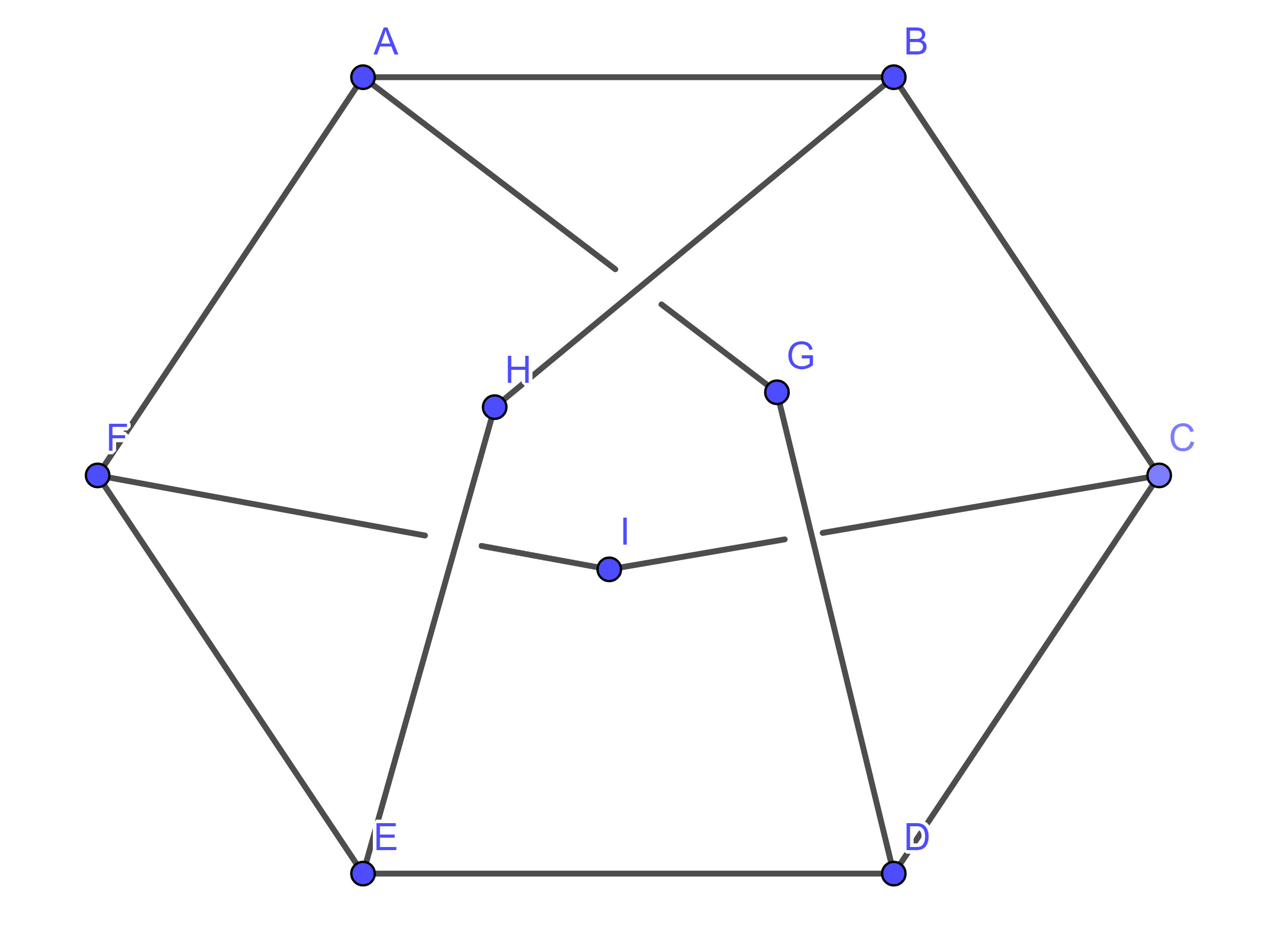} \\
\caption{Example of an indecomposable 3DOF linkage.}
\label{fig:Indecomposable}
\end{figure}

v) A natural question arises:
given a planar linkage $\cL$, what is the maximal number $|I|$ of
components $\cL_i$ ($i\in I$) in a decomposition of $\cL$ into 
acyclic semi-rigid connected sums (with any choice of marked 
diagonals, such that the dimension of the configuration space 
of each $\cL_i$ is strictly greater than the number of marked 
diagonals in it)? If this maximal number is 1, then we 
say that $\cL$ is \textbf{\textit{indecomposable}}. Figure 
\ref{fig:Indecomposable} is an example of a 3DOF indecomposable 
linkage (i.e., the configuration space has dimension 3). 

vi) In general, it is clear that the number of components in a 
``valid'' acyclic semi-rigid connected sum (``valid'' means 
that $\dim M_{\cL_i}$ is strictly greater than the number of 
marked diagonals in $\cL_i$) cannot be greater than $(|\cV| - 1)/2$,  
where $|\cV|$ is the number of vertices of the linkage $\cL$, 
because each time one adds a new component to such an acyclic 
semi-rigid connected sum, one needs to add at least 2 new vertices, 
and the first component has at least 3 vertices. An example where 
the maximal number $(|\cV| - 1)/2$ of components is achieved is 
shown in Figure \ref{fig:Max}. The linkage $\cL$ has 7 vertices, 
the configuration space has 4 dimensions, and the following 
decomposition of $\cL$ into 3 components is valid: $\cL_1$ is 
generated by $\{A,B,G\}$,  $\cL_2$ is generated by $\{B,C,F,G\}$,
$\cL_3$ is generated by $\{C,D,E,F\}$, and the only
marked diagonal is $CF$. So, on this 4-dimensional configuration 
space, we have an integrable system of type $(3,1)$, with three 
commuting vector fields and one common first integral. 

\begin{figure}[!ht]
\centering
\includegraphics[width = 0.6\linewidth]{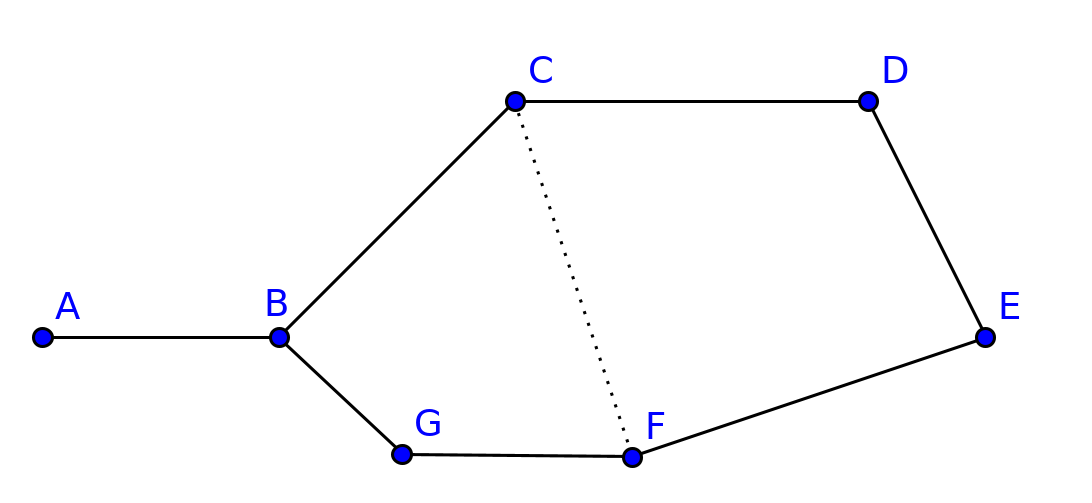} \\
\caption{A linkage which can be decomposed into 3 ``elementary'' components.}
\label{fig:Max}
\end{figure}
}
\end{remark}

\section{Some final remarks and open problems}
\label{section:Final}

\subsection{Bott-Morse theory and graph-manifolds} \hfill

The topology of configuration spaces of planar linkages has been
extensively studied using mostly Morse theory, see, e.g., 
\cite{Farber-Robotics2008,Farber-Homology2007,GiMa-planar-6bar1989,HaKn-3D1997,
HaKn-Cohomology1998,JoSt-Configuration1999,KaTsu-Euler2014,KaMi-moduli1995}
and references therein. However, in general, the square length 
functions that we use in this paper are not Morse functions, but 
rather Bott-Morse functions: their singular points are 
not isolated but form submanifolds and are transversally non-degenerate 
with respect to these submanifolds. So the idea of using Bott-Morse 
functions (instead of just Morse functions) for the study of the 
topology of linkage spaces is very natural.

At least in the case of linkages with 3-dimensional spaces which 
admit integrable systems of type (2,1), this idea works very well: 
the situation is then very similar to Fomenko's topological theory of 
2-degree-of-freedom integrable Hamiltonian systems restricted
to 3-dimensional isoenergy manifolds \cite{Fomenko-2DOF1989}. 
One of the main observations made by Fomenko is that, due to the
Bott-Morse nature of the first integral, the 3-dimensional manifolds 
in question must be \textit{graph-manifolds},  i.e., they can be 
cut into pieces which are direct products of $S^1$ with 2-dimensional 
surfaces with boundary. The class of graph-manifolds has been introduced 
by Waldhausen \cite{Waldhausen-Graphs1967} in 1967  and is a very  
well-understood class of 3-manifolds. (Incidentally, according to the
Thurston-Perelman geometrization program, this class coincides with the
class of 3-manifolds whose Gromov norm vanishes). 

For example, consider the configuration space $M_\ell$ of hexagons 
$A_1A_2A_3A_4A_5A_6$ with fixed edge lengths $|A_iA_{i+1}| = \ell_i$
($i=1,\hdots,6$, $A_7 = A_1$). We assume that the lenghths are generic 
and that  $M_\ell$ is not empty. Denote by $d = A_1A_4$ the length of 
the diagonal $A_1A_4$. Then $d$ is a function on $M_\ell$, which 
is not smooth; however, $\mu = d^2/2$ is smooth and is a Bott-Morse 
function on $M_\ell$, whose singular points form circles.
It follows immediately that $M_\ell$ is a graph-manifold.
To study the topology of $M_\ell$ more carefully, one needs to 
investigate in detail the singularities of $\mu$.
For example, assume (without loss of generality for the topology 
of $M_\ell$) that $\ell_1 \leq \ell_2 \leq \ell_3 \leq
\ell_4 \leq \ell_5  \leq \ell_6$. Then the possible critical 
values of $d$ are: 
\begin{itemize}
\item the maximal value: $\ell_1 + \ell_2 + \ell_3$, 
\item  possible saddle values: $\ell_1 + \ell_2 - \ell_3$ (if this 
value is positive),
$\ell_3 + \ell_2 - \ell_1$, $\ell_3 + \ell_1 - \ell_2$, 
$\ell_4 + \ell_5 - \ell_6$ (if this value is positive),
$\ell_6 + \ell_5 - \ell_4$, $\ell_6 + \ell_4 - \ell_5$;
\item  the minimal value is either 0 (if both 
$\ell_1 + \ell_1 - \ell_3$ and $\ell_4 + \ell_5 - \ell_6$ 
are positive)  or $\max (\ell_3 - \ell_1 - \ell_2, 
\ell_6 - \ell_4 - \ell_5)$ (if this value is positive). 
\end{itemize} 
Then one describes the level sets at and near these critical 
values, the 1-cycles on them, and so on, to recover a long list of all possibilities (see \cite{KaMi-moduli1995}).
This would be an interesting exercise for 
graduate students.

\subsection{Integrable systems on more general  linkage spaces} 
\hfill

In this paper, we used our decomposition and the cross product 
method to construct integrable systems on the configuration 
spaces of planar linkages. It it clear that our method can be 
applied, in a straightforward way, to other classes of linkages, 
including spherical linkages, 3D linkages, and higher-dimensional 
linkages. (See \cite{McCarthy-Designs2011} for a general theory 
of linkage designs, including spherical linkages and spatial 
linkages). In particular, according to our method, the class 
of configuration spaces of 3D linkages admitting interesting 
natural integrable systems is much bigger than the class of 
3D polygon spaces studied by Kapovich-Millson 
\cite{KaMi-3D1996} and others. As far as we know, this bigger 
class has not yet been explored  and it is an interesting 
subject of study, with potential applications in robotics.

\subsection{Singularities of integrable systems on linkage spaces} 
\hfill

Singularities of integrable Hamiltonian systems have been studied 
by many authors. More recently, there has been interest in
formulating a theory of singularities of integrable 
\textit{non-Hamiltonian} systems (see \cite{Zung-NonHamiltonian2016} 
and references therein). A detailed study of singularities of 
concrete integrable systems on configuration spaces of planar 
linkages would help the development of this theory.

\subsection{Commuting flows for 2DOF components} \hfill

In this paper, in the construction of integrable systems 
on configuration spaces, we only used components which are  1DOF 
(i.e., curve-drawing), once the lengths of the marked diagonals 
are fixed. What about 2DOF components? Can they be used effectively? 
In particular, is there any natural way to construct a pair of 
independent commuting vector fields on the configuration spaces 
of pentagons, for example? These configuration spaces are closed 
surfaces and it is known that on such surfaces there exist 
$\bbR^2$-actions with nondegenerate singularities
(see, e.g., \cite{MinhZung-Rn2014}). The question is hence how 
to construct such  $\bbR^2$-actions in a natural way on 
2-dimensional moduli spaces of planar linkages. What
about components with more degrees of freedom (once the
lengths of the marked diagonals are fixed).

\section*{Acknowledgements}

We  thank Marc Troyanov for some interesting discussions on the subject 
of this paper. We also thank B. Servatius, O. Shai, and W. Whiteley 
for their permission to copy a picture from their paper
\cite{Servatius-Assur2010}. Most of this work was carried out while 
the second author visited the first author at Shanghai Jiao Tong University 
for a month in 2019, under its ``High-End Foreign Experts" cooperation program;
he thanks SJTU and the members of the School of Mathematical Sciences, 
especially Xiang Zhang and Jie Hu, for the invitation and their 
warm hospitality.

\end{document}